\documentclass[sn-mathphys,Numbered]{sn-jnl}% Math and Physical Sciences Reference Style
\usepackage{graphicx}%
\usepackage{multirow}%
\usepackage{amsmath,amssymb,amsfonts}%
\usepackage{amsthm}%
\usepackage{mathrsfs}%
\usepackage[title]{appendix}%
\usepackage{xcolor}%
\usepackage{textcomp}%
\usepackage{manyfoot}%
\usepackage{booktabs}%
\usepackage{algorithm}%
\usepackage{algorithmicx}%
\usepackage{algpseudocode}%
\usepackage{listings}%

\begin{document}

\title[Article Title]{Exploring Optimization Techniques for Parameter Estimation in Nonlinear System Modeling}


\author{\fnm{Kaushal} \sur{Kumar}}\email{kaushal.kumar@stud.uni-heidelberg.de}

\affil{\orgdiv{Institute for Applied Mathematics}, \orgname{Heidelberg University}, \orgaddress{\street{Im Neuenheimer Feld}, \city{Heidelberg}, \postcode{69120}, \country{Germany}}}


\abstract{Optimization techniques play a crucial role in estimating parameters and state information for nonlinear systems. However, some critical aspects of these problems have received little attention in previous research. In this paper, we address this gap by exploring optimization techniques for parameter estimation in nonlinear system modeling, with a focus on chaotic dynamical systems. We introduce three optimization methods - a gradient-based iterative algorithm, the Levenberg-Marquardt algorithm, and the Nelder-Mead simplex method - that transfer the complex nonlinear optimization problem into a simpler linear or nonlinear one. We apply these methods to determine the parameters of nonlinear systems, presenting a numerical example to demonstrate their effectiveness. Our results show that the Nelder-Mead simplex method is particularly effective in estimating the parameters of nonlinear systems and has the potential to be a valuable tool in various fields that require nonlinear system modeling. Overall, our study contributes to the understanding and improvement of optimization techniques for parameter estimation in nonlinear system modeling, which has implications for a wide range of applications in science and engineering.}




\maketitle

\section{Introduction}\label{sec1}
Physical models of nonlinear frameworks figured as differential equations or as discrete-time maps regularly have obscure boundaries speaking to our absence of nitty-gritty information or our conjunctures about how the elements work in these frameworks \cite{PhysRevLett.104.060201}.

Nonlinear fitting problems can be solved through iterative methods that minimize a cost function. In system modeling, iterative or recursive algorithms are widely used and have been extensively studied \cite{doi:10.1073/pnas.1517384113,raissi2018multistep,Bock2007,10.1007/978-3-319-23321-5_1}. This paper focuses on the parameter estimation problem of nonlinear systems, which is an important aspect of system identification. System identification is building the mathematical models of systems from the observed data by minimizing a cost function using different optimization methods, just like the gradient methods, Newton methods, or the least square methods \cite{nocedal2006numerical,LI20134278,JOHNSON19921} .\\
To address this problem, we implement optimization methods that simplify the optimization problem by transforming a complex nonlinear system into a linear-in-parameters model or a simpler nonlinear model that is easier to solve. We implement a gradient-based iterative algorithm, the Levenberg-Marquardt algorithm, and the Nelder-Mead simplex method to identify the parameters of a nonlinear function, using an example to demonstrate the effectiveness of these methods. Overall, our study contributes to the understanding and improvement of optimization techniques for parameter estimation in nonlinear system modeling, which has implications for a wide range of applications in science and engineering.\\
The structure of this paper is as follows: In Section 2, we present the objectives of the research problems that we aim to address. Section 3 describes the methods and algorithms used for parameter identification in nonlinear dynamical systems. The results of our experiments are discussed in Sections 4 and 5, where we present the results and conclusions, respectively.

\section{Theory: Background and preliminaries}
Consider a system of ordinary differential equations for state variables $x(t)$ with parameter estimation problems \cite{PhysRevA.45.5524}. The system is described by the following differential equation:
\begin{equation}
\dot{x}(t) = F(x,t,p),
\end{equation}
where $p=(p_{1},p_{2},...,p_{n})$ is the parameter vector. In addition, measurements $\eta_{ij}$ for the state variables or system capacities are provided as:
\begin{equation}
\eta_{ij}= g_{ij}(x(t_{j}),p)+\varepsilon_{ij},
\end{equation}
where $t_{j}$ is the measurement time, $j=1,2,...,k$, and $\varepsilon_{ij}$ is the measurement error. The objective is to find the optimal values of the unknown parameters $p$ such that the model accurately reproduces the observed process. This can be achieved by minimizing an appropriate objective function, which takes into account the measurement errors $\eta_{ij}$. A suitable objective function is the weighted $l_{2}$ norm of the measurement errors, given by:
\begin{equation}
l_{2}(x,p)= \sum_{i,j} \sigma_{ij}^{-2}\varepsilon_{ij}^{2}= \sum_{ij} \sigma_{ij}^{-2}[\eta_{ij}-g_{ij}(x(t_{j}),p)]^{2},
\end{equation}
where $\sigma_{ij}^{2}$ is the variance of the measurement errors, assumed to be independent and Gaussian with zero mean.

To solve this problem, optimization algorithms are used to find the parameter vector $p$ and trajectory $x$ that minimize the objective function. One commonly used algorithm is the gradient-based iterative algorithm and a Newton iterative algorithm, which is an iterative method for solving nonlinear least-squares problems. The algorithm starts with an initial guess for the parameter vector and trajectory and iteratively improves the estimates until a minimum is reached. At each iteration, the algorithm linearizes the model around the current estimate of the trajectory, resulting in a linear least-squares problem. This problem is then solved to obtain an updated estimate of the trajectory and parameter vector. The process continues until the algorithm converges to a minimum, which is a point where the gradient of the objective function is zero. At this point, the estimates for the parameter vector and trajectory are optimal, providing the best fit to the observed data.
\section{ Methods}
\subsection{The gradient-based iterative algorithm}
An algorithm based on gradient estimation can be used to find iterative solutions of $p$ through the gradient search principle. In this algorithm, let $k$ be the iterative variable, while $\hat{p}_{k}$ represents the iterative estimates of $p$ at iteration $k$.
The objective is to create identification techniques for estimating the parameters $p$ by utilizing the available measurement data $(x_{i}; f(x_{i}))$, which is equivalent to minimizing the following cost function
\begin{equation}
    J(p)=\sum_{ij} \sigma_{ij}^{-2}[\eta_{ij}-g_{ij}(x(t_{j}),p)]^{2}.
\end{equation}

and gradient search principle leads to the following gradient-based estimation algorithms
\begin{equation}
    \left[\hat{p}_{k}\right]=\left[\hat{p}_{k-1}\right]-\mu_{k} \left[\frac{\partial J(x,p)}{\partial p} \right]
\end{equation}

where $\mu_{k} > 0$ is the step-size or convergence factor and determined by
\begin{equation}
   \mu_{k}=   argmin\underset{\mu \geq 0}  J(\hat{p_{1},}_{k-1}-\mu \frac{\partial J}{\partial p_{1}}, \hat{p_{2},}_{k-1}-\mu \frac{\partial J}{\partial p_{2}},...) 
\end{equation}
The computation procedure of the gradient-based estimation algorithm for the system described in equation (1) can be summarized as follows:
\begin{enumerate}
\item Collect the measured data $(x_{i},f(x_{i})),$ $i=1,2,...N.$
\item To initialize, let $k=1$, $\hat{p_{0}}$ be arbitrary real numbers, and the pre-set small $\varepsilon$
\item Determine the step-size $\mu_{k}$ by (6).
\item Compute $\hat{p_{k}}$ by (5).
\item If $\sum  \| \hat{p}_{k} -\hat{p}_{k-1} \| > \epsilon$, increase $k$ by $1$ and go to step 3; otherwise, terminate the procedure and obtain the estimate $\hat{p_{k}}$.
\end{enumerate}

The learning rate also needs to be chosen carefully to ensure that the algorithm converges efficiently without oscillating or overshooting the minimum.

\subsection{The Levenberg-Marquardt algorithm}
The Levenberg-Marquardt algorithm \cite{Nocedal20061} iteratively updates the parameter estimates using a combination of the Gauss-Newton method and the steepest descent method, with a damping parameter that adjusts the step size based on the curvature of the cost function. This algorithm is widely used for nonlinear least squares problems and can be applied to a variety of parameter estimation problems in different fields.

\begin{algorithm}[h]
\caption{Levenberg-Marquardt algorithm for parameter estimation problems}
\label{alg:lm}
\begin{algorithmic}[1]
\Require Initial parameter estimates $p_0$, damping parameter $\lambda_0 > 0$, and tolerance $\epsilon > 0$.
\Ensure Best estimate of the parameters $p_0$.
\State Collect the measured data $(x_{i}, f(x_{i}))$ for $i = 1, 2, ..., N$.
\State Define the residual function $e(p)=[f(x_{1})-y_{1},f(x_{2})-y_{2},...,f(x_{N})-y_{N}]$, where $y_{1},y_{2},...,y_{N}$ are the measured response values.
\State Calculate the Jacobian matrix $J(p)= [\frac{\partial f}{\partial p_{1}},\frac{\partial f}{\partial p_{2}},...,\frac{\partial f}{\partial p_{n}}]$, evaluated at the current parameter $p$.
\State Initialize $p = p_0$, $\lambda = \lambda_0$, and $J(p)$.
\While{$\Delta p > \epsilon$}
\State Compute the Gauss-Newton step, $dp_{GN}$, by solving the linear system of equations: $[J(p)^{T} J(p) + \lambda I] dp_{GN} = -J(p)^{T} e(p)$, where $I$ is the identity matrix and $T$ denotes matrix transpose.
\State Calculate the trial parameter estimates $p_{trial} = p + dp_{GN}$.
\State Compute the residual vector $e(p_{trial})$ and evaluate the corresponding cost function $J(p_{trial})$.
\If{$J(p_{trial}) < J(p)$}
\State Update the parameter estimates $p = p_{trial}$, reduce the damping parameter $\lambda$ by a factor of 10, and update $J(p)$.
\Else
\State Increase the damping parameter $\lambda$ by a factor of 10 and update $J(p)$.
\EndIf
\State Calculate $\Delta p = ||p_{trial} - p||$
\EndWhile
\State \Return $p$
\end{algorithmic}
\end{algorithm}

\subsection{The Nelder-Mead simplex method}

The Nelder-Mead simplex method \cite{Nocedal2006} iteratively updates the simplex by reflecting, expanding, and contracting its vertices based on the cost function evaluations. This algorithm is simple to implement and does not require gradient information, making it suitable for nonlinear optimization problems with noisy or discontinuous cost functions. However, it may converge slowly or get stuck in local minima in some cases, so it may be necessary to use multiple starting points or combine them with other optimization techniques.  

\begin{algorithm}[h]
\caption{Nelder-Mead Simplex Method for parameter estimation problems}
\begin{algorithmic}[1]
\Require Initial parameter estimates $p_0$ and simplex $S_0$ with $n+1$ vertices, where $n$ is the number of parameters to be estimated, and covers the parameter space of interest
\Ensure Best vertex $p_0$ as the estimate
\State Evaluate the cost function $J(\pi)$ for each vertex $\pi$ in $S_0$
\State Sort the vertices in the order of increasing cost: $J(p_{0}) \le J(p_{1}) \le \dots \le J(p_{n})$
\While {stopping criterion not satisfied}
\State Compute the centroid of the $n$ best vertices: $c = \frac{1}{n}(p_{0} + p_{1} + \dots + p_{n-1})$
\State Reflect the worst vertex $p_{n}$ about the centroid to obtain the trial point $p_{r} = c + \alpha(c - p_{n})$, where $\alpha > 0$ is a reflection coefficient
\State Evaluate the cost function $J(p_{r})$
\If {$J(p_{r}) < J(p_{0})$}
\State Expand the simplex in the direction of $p_{r}$: $p_{e} = c + \beta(p_{r} - c)$, where $\beta > 1$ is an expansion coefficient. Evaluate the cost function $J(p_{e})$
\ElsIf {$J(p_{r}) \ge J(p_{1})$}
\State Contract the simplex around the best vertex $p_{0}$: $p_{c} = c + \gamma(p_{0} - c)$, where $0 < \gamma < 1$ is a contraction coefficient. Evaluate the cost function $J(p_{c})$
\Else
\State Replace the worst vertex $p_{n}$ with $p_{r}$ and re-sort the vertices
\EndIf
\EndWhile
\end{algorithmic}
\end{algorithm}

\subsection{Accuracy measures}
The Root-Mean-Square Error (RMSE) is commonly used to evaluate the accuracy of a model's predictions. It assesses the variance or residuals between estimated and true values and is utilized to compare the estimated errors of various models for a given model. The formula for calculating RMSE is as follows:
\begin{equation}
    RMSE = \sqrt{\frac{1}{N} \sum_{i=1}^{N} (y^{True} - y^{Estimated})^{2}}
\end{equation}
Where $N$ represents the total number of observations, $y^{True}$ represents the actual value, and $y^{Estimated}$ represents the estimated value.

\section{Results}\label{sec2}
In this section, we conduct numerical tests of the methods introduced above. While the examples we use are not representative of the full range of problems that can be addressed by these methods, they do allow us to evaluate the numerical properties of the algorithms. Specifically, we can assess the stability, reliability, efficiency, and accuracy of the methods, and examine their broader applicability.

We also consider nonlinear ordinary differential equations that can exhibit chaotic solutions for specific parameter values and initial conditions  \cite{strogatz2019nonlinear}  \cite{HIRSCH2013139}. These systems can be challenging to analyze, and testing our methods on such systems provides a rigorous evaluation of their ability to handle complex and nonlinear dynamics. Overall, our numerical tests aim to demonstrate the effectiveness of these methods in solving a range of problems in parameter estimation for nonlinear systems. In this study, we generated simulated data by integrating an ODE and added zero-mean Gaussian noise with a variance of 0.1 to create a realistic dataset. We then employed various optimization algorithms to estimate the system's parameters and successfully captured the underlying dynamics.

\subsection{Test Problem 1.}
\textbf{Lotka Volterra:-}
The exhibition of an ecological system consisting of one predator and one prey may be described by the subsequent model of Lotka and Volterra \cite{ARDITI1989311}:
\begin{eqnarray}
\frac{dx_{1}}{dt} &= &p_{1}x_{1}-p_{2}x_{1}x_{2}\nonumber\\
\frac{dx_{2}}{dt} &=& -p_{3}x_{2}+p_{4}x_{1}x_{2}
\end{eqnarray}

where $x_{1}$, and $x_{2}$ are state variables, and $p_{1},p_{2},p_{3}$, and $p_{4}$ are the system parameters.

\begin{figure}[h]%
\centering
\includegraphics[width=0.99\textwidth]{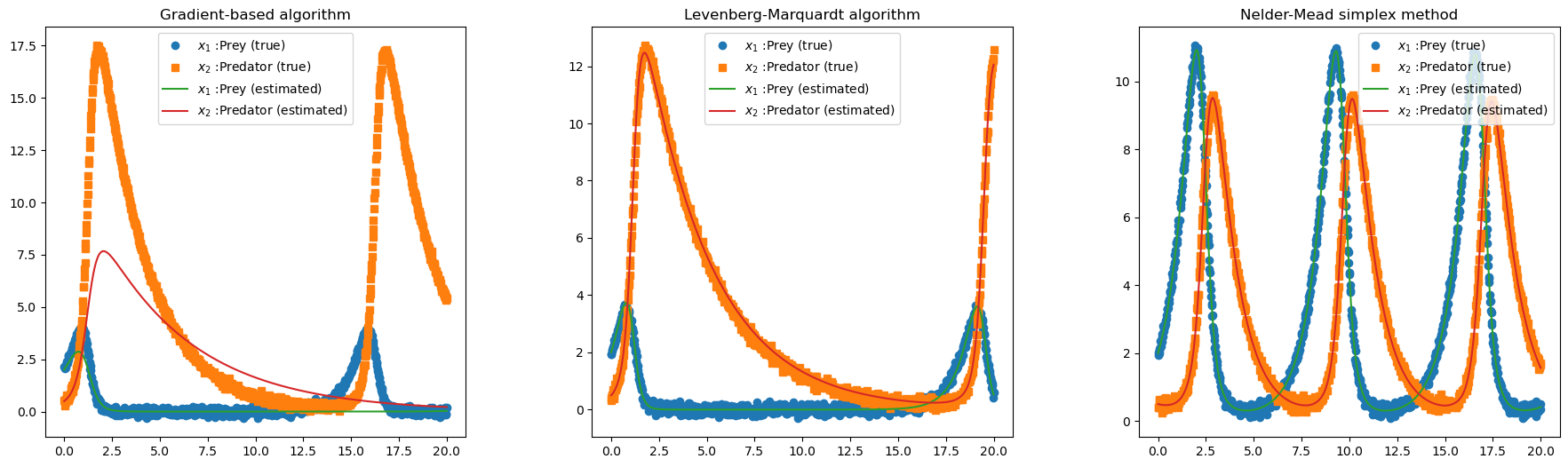}
\caption{The exact trajectories of the Lotka Volterra systems is compared to the
corresponding trajectories of the learned dynamics. Solid blue lines represent the exact
dynamics while the red solid lines demonstrate the learned dynamics}\label{fig1}
\end{figure}

\begin{figure}[h]%
\centering
\includegraphics[width=0.99\textwidth]{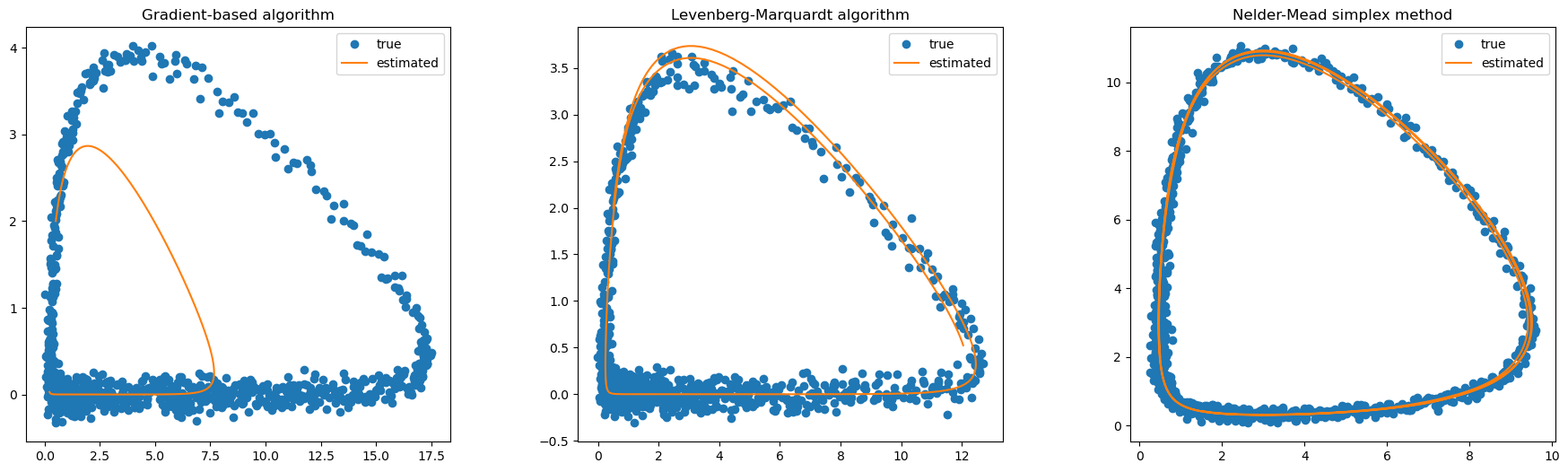}
\caption{The exact phase portrait of the Lotka Volterra systems, where dotted points  are
compared to the corresponding phase portrait of the learned dynamics with all the methods}\label{fig1}
\end{figure}
Based on the initial condition of $[x_{1_{0}} \hspace{0.1cm} x_{2_{0}}]^{T} = [2.0 \hspace{0.2cm} 0.5]^{T}$, data were collected from $t = 0$ to $t = 20$ using a time-step size of $\delta t = 0.01$. In Figure 1, where blue dots and red squares present the collected data, which were used to identify the nonlinear dynamics of the system using different optimization methods. Specifically, iterative gradient-based, Levenberg-Marquardt, and Nelder-Mead simplex methods were employed for parameter estimation. The identified system was then solved using the same initial condition, and the accuracy of the representation of the nonlinear dynamic was assessed using qualitative analysis.\\

\begin{table}[h]
\begin{minipage}{170pt}
\caption{Parameter Identification for Lotka and Volterra}\label{tab1}%
\begin{tabular}{@{}lllll@{}}
\toprule
Parameter & True Value  & Gradient-based  & Levenberg-Marquardt & Nedler-Mead\\
\midrule
$p_{1}$   & 1.2  & 1.0125  & 1.5393 & 1.2014 \\
$p_{2}$   & 0.3  & 0.5218  & 0.4977 & 0.2993\\
$p_{3}$   & 0.4  & 0.2043   & 0.2625 & 0.4004\\
$p_{4}$   & 0.9  & 0.8010   & 0.8555  & 0.8987\\

\botrule
\end{tabular}
\end{minipage}
\end{table}

From the comparison of the true and estimated trajectories of the system and the resulting phase portraits presented in Figure 2, we observed that the Nelder-Mead Simplex algorithm was able to accurately capture the dynamic evolution of the system. The Estimated parameter with Nelder-Mead was also found to be very close to the true parameter as shown in Table 1. \\

\begin{table}[h]
\begin{minipage}{165pt}
\caption{RMSE for Lotka and Volterra with different Algorithms}\label{tab1}%
\begin{tabular}{@{}llll@{}}
\toprule
Methods & Gradient-based & Levenberg-Marquardt & Nedler-Mead\\
\midrule
RMSE   & 4.1373   & 0.2582  & 0.0995  \\
\botrule
\end{tabular}
\end{minipage}
\end{table}

Table 2 presents the RMSE values obtained for the Lotka-Volterra system using different optimization algorithms. The Nelder-Mead simplex method achieved the lowest RMSE value of 0.0995, indicating that it provided the best fit to the true data. In contrast, the iterative gradient-based and Levenberg-Marquardt methods resulted in higher RMSE values of 4.1373 and 0.2582, respectively.

\subsection{Test problem 2.} 
\textbf{Van der Pol oscillator:-}
It is a nonlinear second-order differential equation that exhibits limit cycle behavior \cite{Kanamaru:2007}. The following equation gives it:
\begin{eqnarray}
    x^{''} - \mu (1 - x^{2})x^{'} + x = 0
\end{eqnarray}

where x is the displacement of the oscillator from its equilibrium position, $\mu$ is a parameter that governs the strength of the nonlinear damping, and the prime denotes differentiation with respect to time. Equation (9) can be written as the systems of first-order equations:

\begin{eqnarray}
\frac{dx_{1}}{dt} &= &x_{2}\nonumber\\
\frac{dx_{2}}{dt} &=& \mu (1 - x_{1}^{2})x_{2} - x_{1}
\end{eqnarray}

The Van der Pol oscillator is a canonical example of self-sustained oscillations. 

\begin{figure}[h]%
\centering
\includegraphics[width=0.99\textwidth]{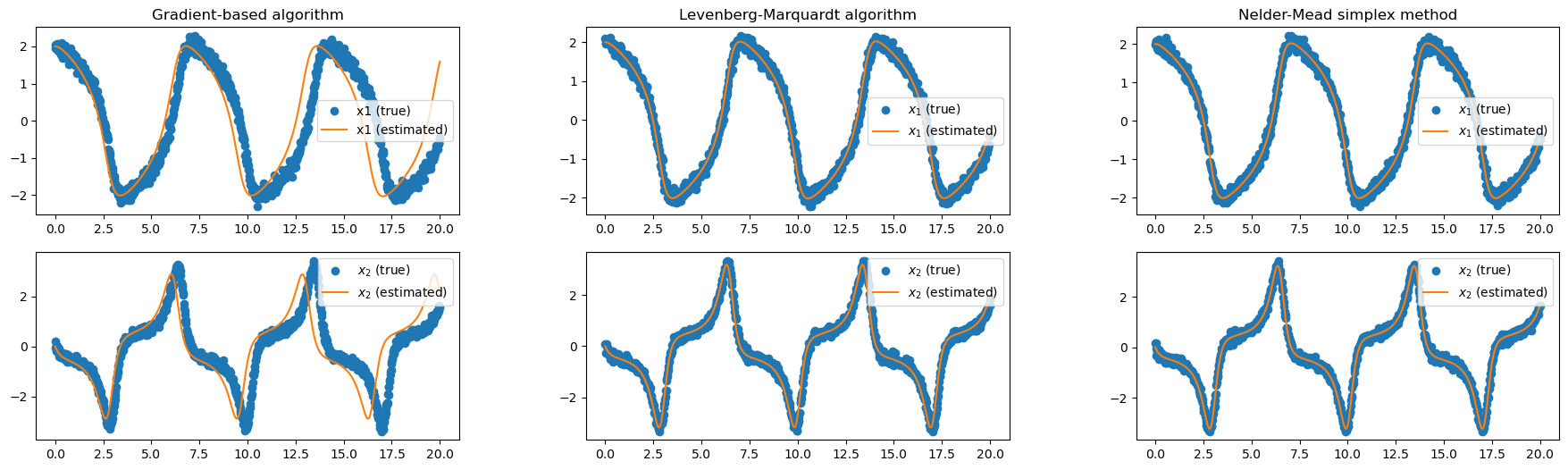}
\caption{The exact trajectories of the Van der Pol oscillator is compared to the
corresponding trajectories of the learned dynamics. Solid blue lines represent the exact
dynamics while the red solid lines demonstrate the learned dynamics}\label{fig1}
\end{figure}

\begin{figure}[h]%
\centering
\includegraphics[width=0.99\textwidth]{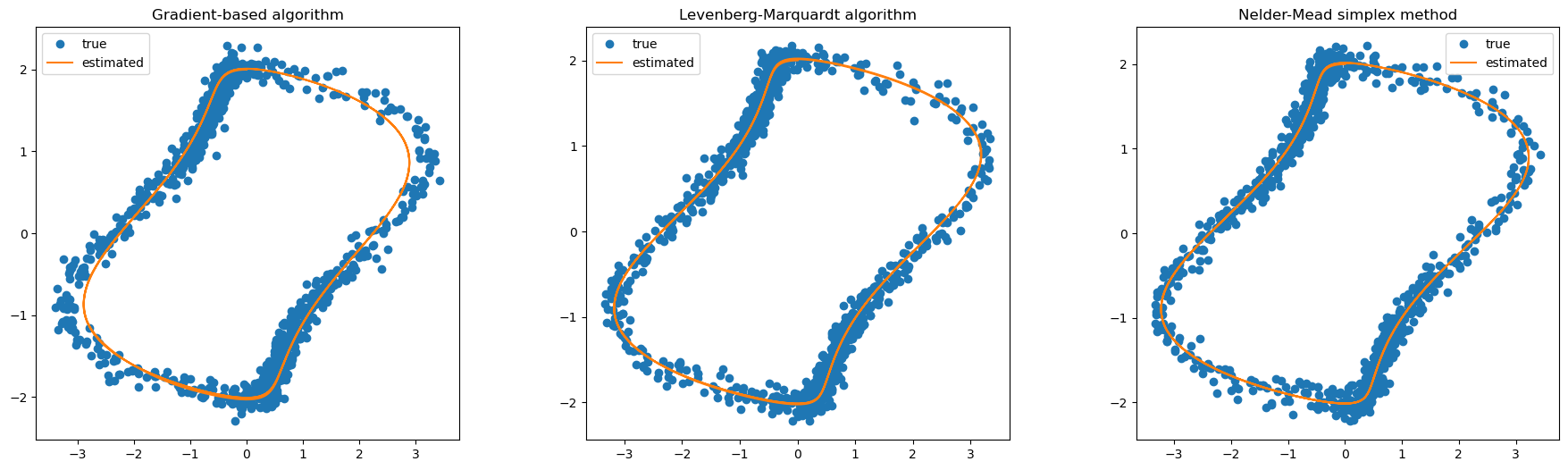}
\caption{The exact phase portrait of the Van der Pol Oscillator, where dotted points  are
compared to the corresponding phase portrait of the learned dynamics with all the methods}\label{fig1}
\end{figure}

\begin{table}[h]
\begin{minipage}{190pt}
\caption{Parameter Identification for Van der Pol oscillator}\label{tab1}%
\begin{tabular}{@{}lllll@{}}
\toprule
Parameter & True Value & Gradient-based & Levenberg-Marquardt & Nedler-Mead\\
\midrule
$\mu$    & 1.5 & 1.2018 & 1.4611  & 1.5007  \\
\botrule
\end{tabular}
\end{minipage}
\end{table}

The initial condition of $[x_{1_{0}} \hspace{0.1cm} x_{2_{0}}]^{T} = [2.0 \hspace{0.2cm} 0.0]^{T}$ was used to collect data from $t = 0$ to $t = 20$ using a time-step size of $\delta t = 0.01$. The collected data, represented by blue dots in Figure 3, were used to identify the nonlinear dynamics of the system using different optimization methods. Specifically, iterative gradient-based, Levenberg-Marquardt, and Nelder-Mead simplex methods were employed for parameter estimation.

The identified system was then solved using the same initial condition, and the accuracy of the representation of the nonlinear dynamic was assessed using qualitative analysis. From the comparison of the true and estimated trajectories of the system and the resulting phase portraits presented in Figure 4, we observed that the Nelder-Mead Simplex algorithm was able to accurately capture the dynamic evolution of the system. The Estimated parameter with Nelder-Mead was found to be very close to the true parameter, as shown in Table 3. \\

\begin{table}[h]
\begin{minipage}{184pt}
\caption{RMSE for Van der Pol Oscillator with different Algorithms}\label{tab1}%
\begin{tabular}{@{}llll@{}}
\toprule
Methods & Gradient-based & Levenberg-Marquardt & Nedler-Mead\\
\midrule
RMSE   & 0.8799  & 0.1409  & 0.1023  \\
\botrule
\end{tabular}
\end{minipage}
\end{table}

Table 4 presents the RMSE values obtained for the Van der Pol oscillator using different optimization algorithms. The Nelder-Mead simplex method achieved the lowest RMSE value of 0.1023, indicating that it provided the best fit to the true data. In contrast, the iterative gradient-based and Levenberg-Marquardt methods resulted in higher RMSE values of 0.8799 and 0.1409, respectively.

\subsection{Test problem 3.}
\textbf{Rossler Systems:-}
The Rössler system is a set of three coupled nonlinear ordinary differential equations that exhibit chaotic behavior. It is named after the German biochemist and mathematician Otto Rössler  \cite{Letellier:2006},\cite{1976ZNatA311664R},\cite{R}, who first introduced it in 1976. The equations are given by:

\begin{eqnarray}
\frac{dx_{1}}{dt} &=& -x_{2}-x_{3}\nonumber\\
\frac{dx_{2}}{dt} &=& x_{1}+p_{1} x_{2}\nonumber\\
\frac{dx_{3}}{dt} &=& p_{2}+x_{3}(x_{1}-p_{3}).
\end{eqnarray}

where $x_{1}, x_{2}$, and $x_{3}$ are the state variables, and $p_{1}, p_{2}$, and $p_{3}$ are the system parameters.\\
The Rössler system is a well-known example of a chaotic system and has been extensively studied in various fields, including physics, engineering, and biology. Its chaotic behavior arises due to the nonlinearity and feedback in the system, which leads to the emergence of a strange attractor. 
\begin{figure}[h]%
\centering
\includegraphics[width=0.99\textwidth]{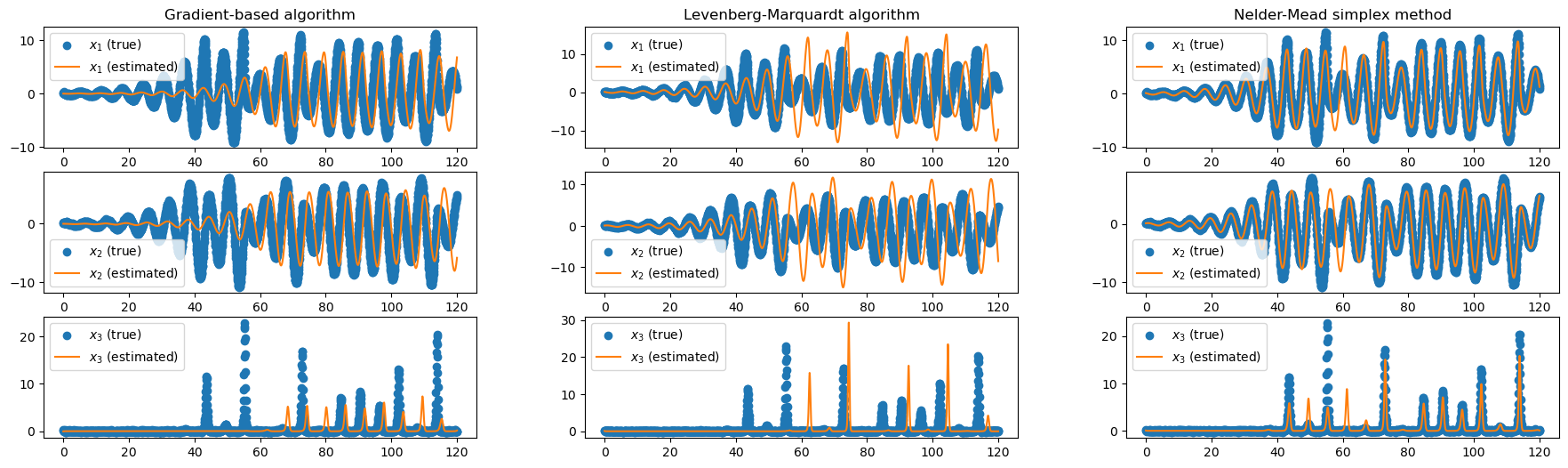}
\caption{The exact trajectories of the Rossler systems is compared to the
corresponding trajectories of the learned dynamics. Solid blue lines represent the exact
dynamics while the red solid lines demonstrate the learned dynamics}\label{fig1}
\end{figure}

\begin{figure}[h]%
\centering
\includegraphics[width=0.99\textwidth]{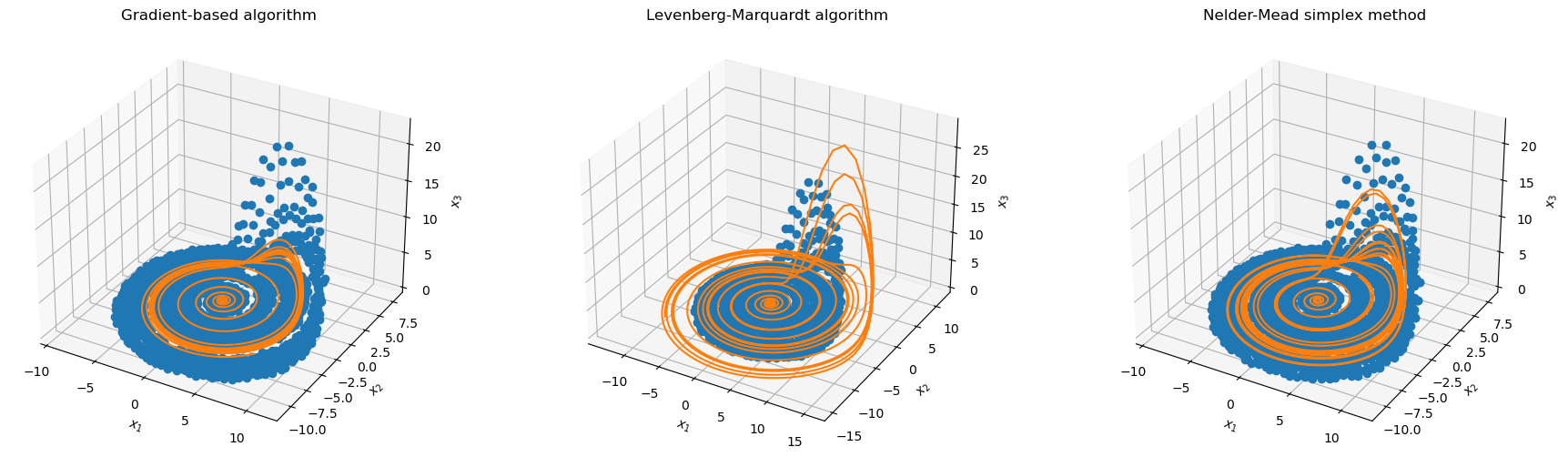}
\caption{The exact phase portrait of the Rossler system dotted points  is
compared to the corresponding phase portrait of the learned dynamics with all the methods}\label{fig1}
\end{figure}

In this study, we aimed to identify the parameters of the Rössler system by collecting data from the system with an initial condition of $[x_{1_{0}} \hspace{0.1cm} x_{2_{0}} \hspace{0.1cm} x_{3_{0}}]^{T} = [0.1 \hspace{0.1cm} 0.1 \hspace{0.1cm} 0.1]^{T}$. We collected data from $t = 0$ to $t = 120$ with a time-step size of $\Delta t = 0.01$, as represented by the blue dots in Figures 5 and 6.\\

\begin{table}[h]
\begin{center}
\begin{minipage}{174pt}
\caption{Parameter Identification for Rössler system}\label{tab1}%
\begin{tabular}{@{}lllll@{}}
\toprule
Parameter &  True Value & Gradient-based & Levenberg-Marquardt & Nedler-Mead\\
\midrule
$p_{1}$  & 0.2 & 0.1711 & 0.1499 & 0.1913\\
$p_{2}$  & 0.2 & 0.1501 & 0.1601 & 0.1918\\
$p_{3}$  & 5.7 & 4.418 & 8.8767 & 4.9344 \\
\botrule
\end{tabular}
\end{minipage}
\end{center}
\end{table}

We employed three different optimization algorithms, namely iterative gradient-based, Levenberg-Marquardt, and Nelder-Mead simplex methods, to estimate the parameters of the Rössler system. Table 5 shows the estimated parameters obtained from the optimization algorithms, with the Nelder-Mead simplex algorithm achieving the closest parameter estimates to the true parameters. We solved the Rössler system using the identified parameters and compared the resulting trajectories and phase portraits with the true system to assess the accuracy of the parameter estimation. Figure 5 and Figure 6 demonstrate the comparison of the true and estimated trajectories of the system and the resulting phase portraits, showing that the Nelder-Mead simplex algorithm was able to accurately capture the dynamic evolution of the system.\\

\begin{table}[h]
\caption{RMSE for Rossler systems with different Algorithms}\label{tab1}%
\begin{tabular}{@{}llll@{}}
\toprule
Methods & Gradient-based & Levenberg-Marquardt & Nedler-Mead\\
\midrule
RMSE   & 4.4744  & 6.9839 & 1.3247  \\
\botrule
\end{tabular}
\end{table}

Furthermore, Table 6 shows the RMSE values obtained for the Rossler system using different optimization algorithms. The Nelder-Mead simplex method resulted in the lowest RMSE value of 1.3247, indicating that it provided the best fit to the true data. On the other hand, the iterative gradient-based and Levenberg-Marquardt methods resulted in higher RMSE values of 4.4744 and 6.9839, respectively.\\

To examine the influence of noisy derivatives in a controlled environment, we introduced zero-mean Gaussian measurement noise with a variance of $\eta$ to the precise derivatives. In Figure 7, we showcase the trajectories of the Rossler system from time t = 0 to t = 120. The true dynamics are represented in red, while the identified systems derived from estimated parameters are illustrated in blue. We assessed the performance of the identified systems under varying levels of additive noise, ranging from $\eta$ = 0.01 to $\eta$ = 10. We observed that the estimated parameters effectively captured the system's dynamics, even under different levels of noise.

\begin{figure}[h]%
\centering
\includegraphics[width=0.99\textwidth]{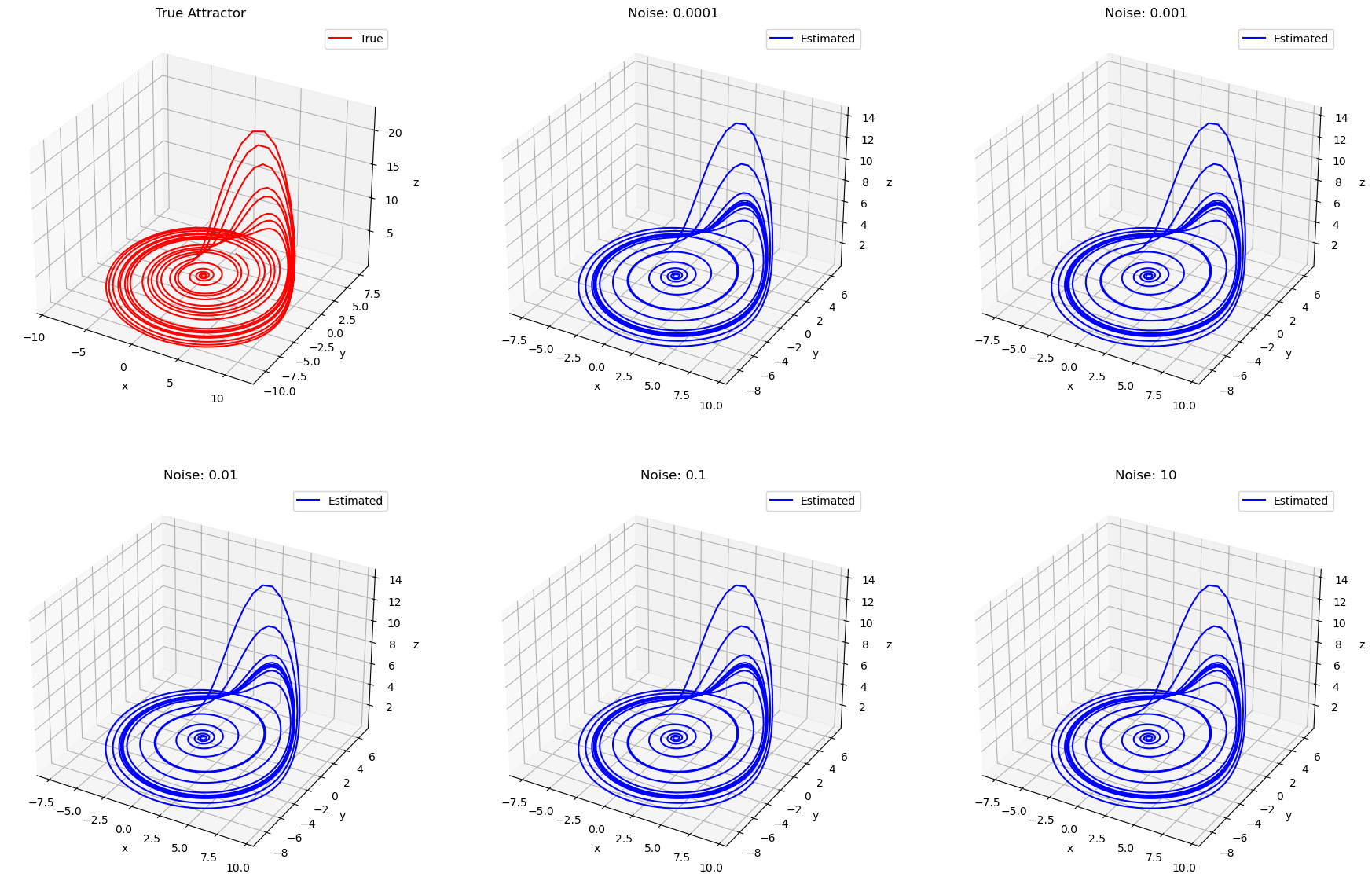}
\caption{In this figure, we present the trajectories of the Rossler system for t = 0 to t = 120. The true dynamics are depicted in red, while the identified systems obtained from estimated parameters are displayed in blue. The performance of the identified systems is evaluated under different levels of additive noise.}\label{fig1}
\end{figure}

\newpage

\section{ Conclusions}

In this study, we investigated the effectiveness of iterative gradient-based, Levenberg-Marquardt, and Nedler-Mead simplex optimization methods for parameter identification in nonlinear dynamical systems. Specifically, we applied these methods to the Rossler system, Van der Pol oscillator, and Lotka-Volterra predator-prey model, and assessed their performance using the root mean square error (RMSE) metric.

Our results show that the Nelder-Mead simplex algorithm consistently outperformed the other two methods in terms of accuracy and efficiency, producing lower RMSE values and accurately capturing the dynamic evolution of the systems. Additionally, our findings indicate that the choice of an optimization algorithm can have a significant impact on the accuracy of parameter estimation and the quality of the resulting system identification.

Overall, our study highlights the importance of carefully selecting an appropriate optimization algorithm for parameter identification in nonlinear dynamical systems. By choosing the right method, researchers can improve the accuracy of their models, better understand the underlying dynamics of complex systems, and make more reliable predictions about their behavior. We hope that our findings will contribute to the development of better methods for system identification and inspire further research in this important area of study.\\

\textbf{Author contribution:} K.K was responsible for the conception and design of the study, as well as overseeing data collection, development of models, analysis, and drafting of the article.\\

\textbf{Code availability:} The codes used to generate the results in this paper are publicly available on GitHub \footnote{\url{https://github.com/kaushalkumarsimmons/Parameter_Estimation1}} following the acceptance of the manuscript. If there is a need for additional information or clarification on the results presented, the corresponding author can be contacted via email, and we will do our best to provide a timely response.\\

\textbf{Funding:} No funding was received for conducting this study.

\section*{Declarations}
\textbf{Conflict of interest: } The authors declare no competing interests.\\


\bibliography{sn-bibliography}% common bib file

\end{document}